\documentclass[12pt,oneside,openany,article]{memoir}
\usepackage{cmap}

\nouppercaseheads
\usepackage[amsthm,thmmarks]{ntheorem}
\usepackage[noTeX]{mmap}
\usepackage{etex}

\usepackage[russian, english]{babel}
\usepackage[leqno]{mathtools}

\usepackage{mathrsfs}
\usepackage{upgreek}
\usepackage[compress,square,comma,numbers]{natbib}
\usepackage{hypernat}

\usepackage{eso-pic}
\usepackage{sfmath}
\usepackage{tikz}

\usepackage{todonotes}

\usepackage{subfig}

\usepackage{microtype}

\usepackage{array}

\setsecheadstyle{\large\bfseries\memRTLraggedright}
\setsubsecheadstyle{\bfseries\memRTLraggedright}

\usepackage{hyperxmp}
\usepackage{xr-hyper}
\usepackage{nameref}
\usepackage[pdftex,bookmarks,pdfnewwindow,plainpages=false,unicode]{hyperref}

\usepackage{bookmark, url}

\usepackage{enumitem}
\usepackage{amssymb}

\newenvironment{claim}[1][{\textup{(\theequation)}}]{\refstepcounter{equation}\vglue10pt
\begin{trivlist}
\item[{\hskip\labelsep#1}]}{\vglue10pt\end{trivlist}}

\hypersetup{
colorlinks=true,
linkcolor=black,
citecolor=black,
urlcolor=blue,
pdfauthor={Victor Ivrii},
pdftitle={Complete Semiclassical Spectral Asymptotics for Periodic and Almost Periodic Perturbations of Constant Operators},
pdfsubject={Sharp Spectral Asymptotics},
pdfkeywords={Microlocal Analysis, sharp  spectral asymptotics, integrated density of states, periodic and almost periodic operators, Diophantine conditions},
bookmarksdepth={4}
}

\theoremstyle{plain}
\newtheorem{theorem}{Theorem}[chapter]

\newtheorem{proposition}[theorem]{Proposition}
\newtheorem{corollary}[theorem]{Corollary}

\theoremstyle{definition}

\theoremstyle{remark}
\newtheorem{remark}[theorem]{Remark}

\newtheorem{condition}{Condition}

\makeatletter
\newtheoremstyle{plainfoot}%
 {\item[\hskip\labelsep \theorem@headerfont ##1\ ##2\,\footnotemark\theorem@separator]}%
 {\item[\hskip\labelsep \theorem@headerfont ##1\ ##2\ (##3)\, \footnotemark\theorem@separator]}
\makeatother

\theoremstyle{plainfoot}
\theoremseparator{.}
\newtheorem{theorem-foot}[theorem]{Theorem}
\newtheorem{lemma-foot}[theorem]{Lemma}
\newtheorem{proposition-foot}[theorem]{Proposition}
\newtheorem{corollary-foot}[theorem]{Corollary}
\newtheorem{conjecture-foot}[theorem]{Conjecture}
\newtheorem{condition-foot}[condition]{Condition}

\theoremstyle{plainfoot}
\theoremseparator{.}
\theorembodyfont{}
\newtheorem{definition-foot}[theorem]{Definition}
\newtheorem{Problem-foot}[theorem]{Problem}

\theoremstyle{plainfoot}
\theoremseparator{.}
\theorembodyfont{}
\theoremheaderfont{\slshape}
\newtheorem{remark-foot}[theorem]{Remark}     
\newtheorem{example-foot}[theorem]{Example}
\newtheorem{problem-foot}[theorem]{Problem}

\numberwithin{equation}{chapter}

\newcounter{note}

\newcommand{\B}{\mathsf{B}}

\newcommand{\bR}{\mathbb{R}}
\newcommand{\bC}{\mathbb{C}}
\newcommand{\bZ}{\mathbb{Z}}

\newcommand{\cN}{\mathcal{N}}

\newcommand{\cV}{\mathcal{V}}
\newcommand{\cZ}{\mathcal{Z}}

\newcommand\sC{\mathscr{C}}
\newcommand\sL{\mathscr{L}}
\newcommand\sH{\mathscr{H}}

\newcommand{\Ad}{\operatorname{Ad}}
\newcommand{\diag}{\operatorname{diag}}
\newcommand{\dist}{\operatorname{dist}}
\newcommand{\diam}{\operatorname{diam}}

\newcommand{\mes}{\operatorname{mes}}

\newcommand{\vol}{\operatorname{vol}}
\newcommand{\supp}{\operatorname{supp}}

\newcommand{\Span}{\operatorname{span}}

\newcommand{\M}{\mathsf{M}}
\newcommand{\N}{\mathsf{N}}

\newcommand {\fU}{\mathfrak{U}}
\newcommand {\fV}{\mathfrak{V}}
\newcommand {\fW}{\mathfrak{W}}

\newcommand{\3}{{|\!|\!|}}

\addtopsmarks{headings}{}{%
\createmark{chapter}{right}{shownumber}{}{. \ }
}
\title{Complete Semiclassical Spectral Asymptotics for Periodic and Almost Periodic Perturbations of Constant Operators\thanks{\emph{2010 Mathematics Subject Classification}: 35P20.}\thanks{\emph{Key words and phrases}: Microlocal Analysis, sharp  spectral asymptotics, integrated density of states, periodic and almost periodic operators, Diophantine conditions.}
}

\author{Victor Ivrii\thanks{This research was supported in part by National Science and Engineering Research Council (Canada) Discovery Grant RGPIN 13827}}

\begin{document}

\maketitle

\begin{abstract}
Under certain assumptions we derive a complete semiclassical asymptotics of the spectral function 
$e_{h,\varepsilon}(x,x,\lambda)$ for a scalar operator
\begin{equation*}
A_\varepsilon (x,hD)= A^0(hD) + \varepsilon B(x,hD),
\end{equation*}
where $A^0$ is an elliptic operator and $B(x,hD)$ is a periodic or almost periodic perturbation.

In particular, a complete semiclassical asymptotics of  the integrated density of states also holds. Further, we consider generalizations.
\end{abstract}


\chapter{Introduction}
\label{sect-1}

\section{Preliminary remarks}
\label{sect-1.1}

This work is inspired by several remarkable papers of L.~Parnovski and R.~Shterenberg \cite{ParSht1, ParSht2, ParSht3}, S.~Morozov, L.~Parnovski and R.~Shterenberg \cite{MorParSht} and earlier papers by A.~Sobolev \cite{Sob1, Sob2}. I wanted to understand the approach of the authors and, combining their ideas with my own approach, generalize their results.

In these papers the complete asymptotic expansion of the integrated density of states $\N(\lambda)$ for operators $\Delta +V$  was derived as $\lambda\to +\infty$; here  $\Delta$ is a positive Laplacian and $V$ is a periodic or almost periodic potential (satisfying certain conditions). In \cite{MorParSht} more general operators were considered.

Further, in \cite{ParSht3} the complete asymptotic expansion of $e(x,x,\lambda)$ was derived, where $e(x,y,\lambda)$
is the Schwartz kernel of the spectral projector.

I borrowed from these papers Conditions~\ref{cond-A}--\ref{cond-D} and the \emph{special gauge transformation\/} and added the \emph{hyperbolic operator method\/} (actually  non-stationary semiclassical Schr\"odinger operator method--\cite{IvrMB}) and extremely long propagation of singularities. I believe that this is a simpler and more powerful approach. Also, in contrast to those papers I consider more general semiclassical asymptotics.

Consider a scalar self-adjoint $h$-pseudo-differential operator $A(x,hD)$ in $\bR^d$ with the Weyl symbol $A(x,\xi)$, such that
\begin{gather}
|D^\alpha _x D^\beta_\xi A(x,\xi)|\le c_{\alpha\beta}(|\xi|+1)^m \qquad \forall \alpha,\beta, \ \forall x,\xi
\label{eqn-1.1}\\
\shortintertext{and}
A(x,\xi)\ge c_0 |\xi|^m - C_0 \qquad \qquad \ \forall x,\xi.
\label{eqn-1.2}
\end{gather}

Then it is semibounded from below. Let $e_h(x,y,\lambda)$ be the Schwartz kernel of its spectral projector
$E(\lambda)=\uptheta (\lambda -A)$. We are interested in the semiclassical asymptotics of  $e_h(x,x,\lambda)$ and
\begin{equation}
\N_h(\lambda)= \M [e(x,x,\lambda)]\ \coloneqq \lim _{\ell \to \infty} (\mes (\ell X))^{-1} \int_{\ell X }
e(x,x,\lambda)\,dx,
\label{eqn-1.3}
\end{equation}
where $0\in X$ is an open domain in $\bR^d$.  The latter expression in the cases we are interested in does not depend on $X$ and is called \emph{Integrated Density of States\/}.

It is well-known that under  \emph{$\xi$-microhyperbolicity condition on the energy level $\lambda$\/}
\begin{equation}
|A(x,\xi,h) - \lambda| +|\nabla_\xi A(x,\xi,h)|\ge \epsilon_0
\label{eqn-1.4}
\end{equation}
the following asymptotics holds
\begin{gather}
e_h(x,x,\lambda)= \kappa_0(x,\lambda)h^{-d}+O(h^{1-d})\qquad\text{as\ \ } h\to +0,
\label{eqn-1.5}\\
\intertext{and therefore}
\N_h(\lambda)=\bar{\kappa}_0(\lambda)h^{-d}+O(h^{1-d}),
\label{eqn-1.6}\\
\shortintertext{where here and below}
\bar{\kappa}_n(\lambda)=\M [\kappa_n(x,\lambda)].
 \label{eqn-1.7}
\end{gather}
For generalization to matrix operators and degenerate scalar operators see Chapters~\ref{monsterbook-sect-4} and~\ref{monsterbook-sect-5} respectively of \cite{IvrMB}. Also there one can find slightly sharper two-term asymptotics under non-periodicity conditions.

Also it is known (see Chapter~\ref{monsterbook-sect-4}   of \cite{IvrMB}) that under microhyperbolicity condition (\ref{eqn-1.4}) for $|\tau-\lambda|<\epsilon$ the following complete asymptotics holds:
\begin{equation}
F_{t\to h^{-1}\tau} \bigl(\bar{\chi}_T(t) \bigl(Q_{2x} u_h(x,y, t)\, ^tQ_{1y}\bigr)|_{y=x}\bigr) \,\sim \,
\sum_{n\ge 0} \kappa_{n,Q_1,Q_2}'(x,\tau)h^{1-d+n},
\label{eqn-1.8}
\end{equation}
where $u_h(x,y,t)$ is the Schwartz kernel of of the \emph{propagator\/} $e^{ih^{-1}tA}$, $\bar{\chi}\in \sC_0^\infty ([-1,1])$, $\bar{\chi}(t)=1$ at $[-\frac{1}{2},\frac{1}{2}]$, $T\in [h^{1-\delta},T^*]$, $T^*$ is a small constant here and $Q_j=Q_j(x,hD)$ are $h$-pseudo-differential operator; we write operators, acting with respect to  $y$ on Schwartz kernels to the right of it.

Further, it is known that
\begin{gather}
\supp(Q_1)\cap\supp(Q_2)=\emptyset \implies \kappa'_{n,Q_1,Q_2}(x,\tau)=0,
\label{eqn-1.9}\\
\intertext{where $\supp (Q_j)$ is a support of its symbol $Q_j(x,\xi)$ and}
\tau \le \tau^* =\inf _{x,\xi} A(x,\xi) \implies \kappa'_{n,Q_1,Q_2}(x,\tau)=0.
\label{eqn-1.10}\\
\shortintertext{Let}
 \kappa_{n,Q_1,Q_2}(x,\tau)=\int_{-\infty}^\tau \kappa'_{n,Q_1,Q_2}(x,\tau')\,d\tau.
 \label{eqn-1.11}
\end{gather}
In what follows we skip subscripts $Q_j=I$.

\begin{remark}\label{rem-1.1}
This equality (\ref{eqn-1.8}) plus H\"ormander's Tauberian theorem imply the remainder estimates $O(h^{1-d})$ for
$Q_{2x}e_h(x,y,\tau)^tQ_{1y}|_{x=y}$. On the other hand, if we can improve (\ref{eqn-1.8}) by increasing $T^*$, we can improve the remainder estimate to $O(T^{*\,-1}h^{1-d})$\,\footnote{\label{foot-1} Provided $T^*=O(h^{-M})$ for some $M$.}\footnote{\label{foot-2} This plus estimate for $\kappa'_0$ is a major method for obtaining sharp remainder estimates in \cite{IvrMB}.}.
\end{remark}

Observe that for $A=A(hD)$
\begin{gather}
e_h(x,x,\lambda)= \N_h(\lambda) =\kappa_0(\lambda)h^{-d}.
\label{eqn-1.12}\\
\intertext{In this paper we consider}
A(x,hD)=A^0(hD)+\varepsilon B(x,hD),
\label{eqn-1.13}
\end{gather}
where $A^0(\xi)$ satisfies (\ref{eqn-1.1}), (\ref{eqn-1.2}) and (\ref{eqn-1.4}) and $B(x,\xi)$ satisfies (\ref{eqn-1.1}) and   $\varepsilon >0$ is a small parameter.  Later we assume that $B(x,hD)$ is almost periodic and impose other conditions.

First, we claim that for operator (\ref{eqn-1.13}) with $\varepsilon \le \epsilon_0$ the equality (\ref{eqn-1.8}) holds with  $T^* =\epsilon_1 \varepsilon^{-1}$ where $\epsilon_{j}$ are small constants and we assume that  $\varepsilon \ge h^M$ for some $M$. Then  the remainder estimate is $O(\varepsilon h^{1-d})$\,\footnote{\label{foot-3} See Theorem~\ref{thm-2.4}.}.

\section{Main Theorem}
\label{sect-1.2}

Now we consider the main topic of this work where we will use ideas from \cite{ParSht1, ParSht2, ParSht3, MorParSht}: the case  of an almost periodic operator $B(x,hD)$,
\begin{equation}
B(x,\xi)=\sum_{\theta \in \Theta} b_{\theta}(\xi)e^{i\langle \theta,x\rangle}
\label{eqn-1.14}
\end{equation}
with discrete (i.e. without any accumulation points)  \emph{frequency set\/}  $\Theta$.

Operator $B$ is \emph{quasiperiodic\/} if $\Theta$ is a finite set, \emph{periodic\/} if $\Theta$ is a lattice and \emph{almost periodic\/} in the general case.

Our goal is to derive (under certain assumptions) complete semiclassical asymptotics:
\begin{equation}
e_{h,\varepsilon}(x,x,\tau) \sim \sum_{n\ge 0}\kappa_{n,\varepsilon}x(x,\tau)h^{-d+n}.
\label{eqn-1.15}
\end{equation}

First, in addition to microhyperbolicity condition (\ref{eqn-1.4}) we assume that
$\Sigma_\lambda =\{\xi\colon A^0(\xi)=\lambda\}$ is  a \emph{strongly convex surface\/} i.e.
\begin{equation}
\pm \sum_{j,k} A^0_{\xi_j\xi_k}(\xi)\eta_j\eta_k \ge \epsilon|\eta|^2\qquad
 \forall \xi\colon A^0(\xi)=\lambda\ \ \forall \eta\colon \sum_j  A^0_{\xi_j}(\xi)\eta_j=0,
\label{eqn-1.16}
\end{equation}
where the sign depends on the connected component of $\Sigma_\lambda$, containing $\xi$.

Without any loss of generality we assume that \begin{claim}\label{eqn-1.17}
$\Theta$ spans $\bR^d$,  contains $0$ and is symmetric about $0$.
\end{claim}

\begin{condition}\label{cond-A}
For each $\theta_1,\ldots, \theta_d\in \Theta$ \underline{either} $\theta_1,\ldots, \theta_d$ are  linearly independent over $\bR$ \underline{or} they linearly dependent over $\bZ$.
\end{condition}

Assume also that

\begin{condition}\label{cond-B}
For any arbitrarily large $L$ and for any sufficiently large real number $\omega$ there are a finite symmetric about $0$ set
$\Theta'\coloneqq \Theta'_{(L,\omega)}\subset (\Theta \cap \B(0,\omega))$ (with $\B(\xi, r)$ the ball of the radius $r$ and center $\xi$) and a ``cut-off'' coefficients $b'_{\theta}\coloneqq b'_{\theta,(L,\omega)}$, such that
\begin{gather}
B'\coloneqq B'_{( L,\omega)}(x,\xi)\coloneqq \sum_{\theta\in\Theta ' } b'_{\theta} (\xi) e^{i\langle \theta, x\rangle}
\label{eqn-1.18}\\
\shortintertext{satisfies}
\| D_x^\alpha D_\xi^\beta \bigl(B -B' \bigr) \|_{\sL^\infty}\le   \omega^{-L}(|\xi|+1)^m\qquad
\forall \alpha, \beta\colon
 |\alpha| \le L, |\beta|\le L.
 \label{eqn-1.19}
\end{gather}
\end{condition}

\begin{remark}\label{rem-1.2}
\begin{enumerate}[label=(\roman*), wide, labelindent=0pt]

\item\label{rem-1.2-i}
Then
\begin{gather}
|D_\xi^\beta b_{\theta}|= O(|\theta|^{-\infty}(|\xi|+1)^m)\qquad \text{as\ \ }|\theta| \to \infty
\label{eqn-1.20}\\
\shortintertext{and}
|D_\xi ^\beta (b_{\theta} - b'_{\theta})|= O(\omega^{-\infty}(|\xi|+1)^m).
\label{eqn-1.21}
\end{gather}

Indeed, one suffices to observe that $b_\theta(\xi) =\M (B (x,\xi) e^{-i\langle \theta,x\rangle})$ etc.

\item\label{rem-1.2-ii}
On the other hand, under additional assumption
\begin{equation}
\# \{\theta \in \Theta,\, |\theta|\le \omega\} = O(\omega^{p})\qquad \text{as\ \ }\omega \to \infty
\label{eqn-1.22}
\end{equation}
for some $p$, (\ref{eqn-1.20}) implies Condition \ref{cond-B} with $\Theta'_{(L,\omega)}\coloneqq \Theta \cap \B(0,\omega)$. However we will need $\Theta'_{(L,\omega)}$ in the next condition.

\item\label{rem-1.2-iii}
We need only to estimate  the operator norm of the difference between $B(x,hD)$ and $B'(x,hD)$ (from $\sH^m$ to $\sL^2$); therefore for differential operators we can weaken (\ref{eqn-1.19}): if
\begin{equation}
B= \sum_{\mu, \nu: |\alpha|\le m', |\beta|\le m'} D^\alpha b_{\alpha\beta}(x) D^\beta,\qquad b_{\alpha\beta}=b^\dag_{\beta\alpha},
\label{eqn-1.23}
\end{equation}
where we assume that $b_{\alpha\nu}(x)$ and $b'_{\alpha\beta}(x)$ have similar decompositions (\ref{eqn-1.14}) and (\ref{eqn-1.18}) respectively,  then (\ref{eqn-1.19}) should be replaced by
\begin{equation}
\| D_x^\alpha \bigl(b_{\alpha\beta} -b'_{\alpha\beta} \bigr) \|_{\sL^\infty}\le   \omega^{-L} \qquad
\forall \alpha  .
 \label{eqn-1.24}
 \end{equation}

\item\label{rem-1.2-iv}
While Condition \ref{cond-B} is Condition B of \cite{ParSht3}, adopted to our case, Condition \ref{cond-A} and
Conditions \ref{cond-C}, \ref{cond-D} below are borrowed without any modifications (except changing notations).
\end{enumerate}
\end{remark}

The next condition we need to impose is a version of the Diophantine condition on the frequencies of $B$. First, we need some definitions. We fix a natural number $K$ (the choice of $K$ will be determined later by how many terms in the asymptotic decomposition of $e(x,x,\lambda)$ we want to obtain) and consider
$\Theta'_{K}$, which here and below denotes the algebraic sum of $K$ copies of $\Theta'$:
\begin{equation}
\Theta'_K \coloneqq \sum_{1\le i \le K} \Theta.
 \label{eqn-1.25}
 \end{equation}

We say that $\fV$ is a \emph{quasi-lattice subspace} of dimension $q$, if $\fV$ is a linear span of $q$ linear independent vectors $\theta_1,\ldots,\theta_q\in \Theta '_K\setminus 0$. Obviously, the zero space  is a quasi-lattice subspace of dimension $0$ and $\bR^d$ is a quasi-lattice subspace of dimension $d$.

We denote by $\cV_q$ the collection of all quasi-lattice subspaces of dimension $q$ and also
$\cV:=\bigcup_{q\ge 0}\cV_q$.

Consider $\fV,\fU\in\cV$. We say that these subspaces are \emph{strongly distinct}, if neither of them is a
subspace of the other one.  Next, let $\widehat{(\fV,\fU)}\in [0,\pi/2]$  be the angle between them, i.e. the angle between $\fV\ominus\fW$ and $\fU\ominus\fW$, $\fW=\fU\cap\fV$. This angle is positive iff $\fV$ and $\fU$ are strongly distinct.

\begin{condition}\label{cond-C}
For each fixed $L$ and $K$ the sets $\Theta'_{(L,\omega)}$ satisfying (\ref{eqn-1.18}) and (\ref{eqn-1.19}) can be chosen in such a way that for sufficiently large $\omega$ we have
\begin{gather}
s(\omega)=s(\Theta'_K)\coloneqq \inf_{\fV,\fU\in \cV}\sin(\widehat{(\fV,\fU)})\ge \omega^{-1}
\label{eqn-1.26}\\
\shortintertext{and}
r(\omega) \coloneqq \inf_{\theta\in \Theta'_K\setminus 0}|\theta| \ge \omega^{-1},
\label{eqn-1.27}
\end{gather}
where  the implied constant (i.e. how large should $\omega$ be)  depends on $L$ and $K$.
\end{condition}

Let $\fV$ be the span of $\theta_1,\ldots,\theta_{q}\in \Theta'_K\setminus 0$. Then due to Condition~\ref{cond-A} each element of the set
$\Theta'_K\cap\fV$ is a linear combination of $\theta_1,\ldots,\theta_{q}$ with rational coefficients.
Since the set $\Theta'_K\cap\fV$ is finite, this implies that the set
$\Theta'_\infty\cap\fV$ is discrete and is, therefore, a lattice in $\fV$. We denote this lattice by
$\Gamma(\omega;\fV)$.

Our final condition states that this lattice cannot be too dense.

\begin{condition}\label{cond-D}
We can choose $\Theta'_{(L;\omega)}$, satisfying Conditions~\ref{cond-B} and \ref{cond-C}
in such a way that for sufficiently large $\omega$ and for each $\fV\in\cV$, $\fV\ne\bR^d$, we have
\begin{equation}
\vol(\fV/\Gamma(\omega;\fV))\ge\omega^{-1}.
\label{eqn-1.28}
\end{equation}
\end{condition}

\begin{remark}\label{rem-1.3}
See Section 2 of \cite{ParSht3} for discussion of these conditions. In particular, if $\Theta$ is a lattice, then Conditions \ref{cond-A}--\ref{cond-D} are fulfilled. Further, if $\Theta $ is a finite set and Condition \ref{cond-A} is fulfilled, then $\Theta_\infty \coloneqq \bigcup_{K\ge 1} \Theta_K$ is a lattice and Conditions \ref{cond-B}--\ref{cond-D} are fulfilled. Furthermore, the same is true, if $\Theta$ is an arithmetic sum of a finite set and a lattice.
\end{remark}

The main theorem of this paper is
\begin{theorem}\label{thm-1.4}
Let $A$ be a self-adjoint operator \textup{(\ref{eqn-1.13})},  where $A^0$ satisfies
\textup{(\ref{eqn-1.1})}, \textup{(\ref{eqn-1.2})}, \textup{(\ref{eqn-1.4})} and \textup{(\ref{eqn-1.16})} and  $B$ satisfies \textup{(\ref{eqn-1.1})}.

Let Conditions \ref{cond-A}--\ref{cond-D} be fulfilled. Then for $|\tau-\lambda|<\epsilon$, $\varepsilon \le h^{\vartheta}$, $\vartheta  >0$
\begin{equation}
e_{h,\varepsilon}(x,x,\tau  ) \sim \sum_{n\ge 0}\kappa_n(x,\tau ;\varepsilon) h^{-d+n}.
\label{eqn-1.29}
\end{equation}
\end{theorem}

\begin{corollary}\label{cor-1.5}
In the framework of Theorem~\ref{thm-1.4}
\begin{equation}
\N_{h,\varepsilon}(\tau ) \sim \sum_{n\ge 0}\bar{\kappa}_n(\tau ;\varepsilon ) h^{-d+n}.
\label{eqn-1.30}
\end{equation}
\end{corollary}

\section{Plan of the paper}
\label{sect-1.3}

Section~\ref{sect-2} is devoted to the proof of Theorem~\ref{thm-1.4}. In Subsection~\ref{sect-2.1} we make some general remarks, and, in particular, we prove more general albeit far less precise Theorem~\ref{thm-2.4}. Then, in Subsection~\ref{sect-2.2} we describe a gauge transformation.

In Subsection~\ref{sect-2.3} we consider a non-resonant zone and justify such transformation, which reduces operator microlocally to a constant symbol operator $A''(hD, h)$. This allows us to study a propagation of singularities with respect to $\xi$ and prove that the singularities  do not propagate with respect to $\xi$\,\footnote{\label{foot-4} For  time $T^*=h^{-M}$ with arbitrarily large $M$.}. In Subsection~\ref{sect-2.4} we consider a resonant zone and justify such transformation, which reduces operator microlocally to an  operator $A''(x', hD, h)$, where $x'\in \fV$ the corresponding resonant subspace, and prove that the singularities propagate only with respect to $\xi'$. Then the convexity condition implies that the singularities actually do not propagate with respect to $\xi$\,\footref{foot-4}.

In Subsection~\ref{sect-2.5} we consider propagation with respect to $x$ and using the results of Subsections~\ref{sect-2.3} and~\ref{sect-2.4} we prove that the singularities ``propagate away'' and do not return\footref{foot-4}. The we apply Tauberian theorem with $T=T^*$ and prove Theorem~\ref{thm-1.4}.

In Section~\ref{sect-3} we generalize Theorem~\ref{thm-1.4}. First, in Subsection~\ref{sect-3.1} we consider matrix operators with the simple eigenvalues of $A^0(\xi)$.

Then, in Subsection~\ref{sect-3.2} we consider operators  $A^0(hD)+\varepsilon V(x,hD)$ where symbol $V(x,\xi)$ decays as $|x|\to \infty$ and  hybrid operators $A^0(hD)+\varepsilon (B(x,hD)+V(x,hD))$ with almost periodic $B$ and decaying $V$ and show that our methods work for them as well.

Finally, in Subsection~\ref{sect-3.3} we discuss differentiability of our asymptotics with respect to $\tau$.

\chapter{Proof of the Main Theorem}
\label{sect-2}

\section{Preliminary Analysis}
\label{sect-2.1}

\begin{remark}\label{rem-2.1}
\begin{enumerate}[label=(\roman*), wide, labelindent=0pt]
\item\label{rem-2.1-i}
It follows  from Section~\ref{monsterbook-sect-4}  of \cite{IvrMB} that the contribution  of the zone
$\{\xi\colon |A^0(\xi)-\tau|\ge C_0\varepsilon + h^{1-\varsigma}\}$  to the remainder is negligible. Here and below
$\varsigma>0$ is an  arbitrarily small exponent. Namely, let $Q_j=Q_j(hD)$ be operators with the symbols $Q_j(\xi)$, such that
\begin{gather}
\supp(Q_1)\cap\supp(Q_2) \cap \Omega_\tau =\emptyset
\label{eqn-2.1}\\
\shortintertext{with}
\Omega_\tau \coloneqq \{\xi\colon |A^0(\xi)-\tau|\le C_0\varepsilon + h^{1-\varsigma}\}
\label{eqn-2.2}\\
\shortintertext{and satisfying}
|D^\alpha Q_j |\le C_\alpha h^{-(1-\varsigma)|\alpha|}\qquad \forall \alpha.
\label{eqn-2.3}
\shortintertext{Then}
(Q_{2x} e(x,y,\tau) \,^tQ_{1y})|_{y=x} = \kappa_{0,Q_1,Q_2} h^{-d} + O(h^\infty)
\label{eqn-2.4}\\
\shortintertext{with}
\kappa_{0,Q_1,Q_2} =(2\pi)^{-d}\int \uptheta (\tau -A^0(\xi))Q_1(\xi)Q_2(\xi)\,d\xi
\label{eqn-2.5}
\end{gather}
with  $\uptheta (\tau -A^0(\xi))$ equal to either $0$ or $1$ on each connected component of $\Omega_\tau \cap \supp (Q_1)\cap \supp (Q_2))$.

Therefore we  restrict ourself by the analysis in the zone $\Omega_\tau$.
\item\label{rem-2.1-ii}
To upgrade (\ref{eqn-1.8}) with $T=T_*$ (a small constant) to (\ref{eqn-1.8}) with $T=T^*$ it is sufficient to prove that
\begin{equation}
|F_{t\to h^{-1}\tau} \bigl(\chi_T(t) \bigl(Q_{2x} u_h(x,y, t)\,^t\!Q_{1y}\bigr)\bigr|_{y=x}\bigr) |\le C_s h^{-d+s},
\label{eqn-2.6}
\end{equation}
for $|\tau -\lambda|\le \epsilon $, $T\in [T_*, \,T^*]$ and
$\chi\in \sC_0^\infty ([-1,-\frac{1}{2}]\cup[\frac{1}{2},1])$, where $s$ is an arbitrarily large exponent. \end{enumerate}
\end{remark}

In the very general setting for  $|t|\le  h^{-M}$ the propagation speed with respect to $\xi$ does not exceed  $C\varepsilon$. More precisely

\begin{proposition}\label{prop-2.2}
Let $A=A^0+ \varepsilon B$ where $A^0(hD)$ and $B(x,hD)$ are matrix operators satisfying \textup{(\ref{eqn-1.1})}. Let $Q_j(hD)$ be operators with symbols satisfying \textup{(\ref{eqn-2.3})}.  Further, let $\supp (Q_j)\subset \{\xi\colon |\xi|\le c\}$ and
\begin{gather}
\dist (\supp (Q_1),\,\supp (Q_2)) \ge \max (C_0 \varepsilon T,\, h^{1-\varsigma})
\label{eqn-2.7}\\
\intertext{with $T\le h^{-M}$. Then for $|t|\le T$}
\| Q_2 e^{i h^{-1}t A}Q_1\| \le C_{M,s} h^s.
\label{eqn-2.8}
\end{gather}
\end{proposition}

\begin{proof}
One can prove easily by arguments of the proof of Theorem~\ref{monsterbook-thm-2-1-2}  of \cite{IvrMB}, applied to operator  $\varepsilon^{-1}A= \varepsilon^{-1}A^0(hD_x) + B(x,hD)$ and $\phi(\xi, t)$, that the propagation speed with respect $\xi$ does not exceed $C_0$; presence of the term $\varepsilon^{-1}A^0(hD_x)$ does not matter since it disappears in the commutator with $\phi(hD)$. Changing $t\mapsto \varepsilon t$ we conclude that for operator $A$ the propagation speed with respect to $\xi$ does not exceed $C_0\varepsilon$.

We do not need compactness of the domain in the phase space with respect to $x$ since the propagation speed with respect to $x$ does not exceed $C_0$ and we have such compactness implicitly. We leave easy details to the reader.
\end{proof}

\begin{proposition}\label{prop-2.3}
In the framework of Proposition~\ref{prop-2.2} assume that $A^0(hD)$ is microhyperbolic on the energy level $\lambda$\,\footnote{\label{foot-5} For definition for matrix operators see Definition~\ref{monsterbook-def-2-1-1} of \cite{IvrMB}.}.

Then for  $T_*\le T\le T^* =\min(\epsilon_0 \varepsilon^{-1}, h^{-M})$ \textup{(\ref{eqn-2.6})} holds.
\end{proposition}

\begin{proof}
It is sufficient to prove for $\supp(Q_1)$ contained in the small vicinity of some point $\bar{\xi}$. Then due to Proposition~\ref{prop-2.2} $e^{ih^{-1}tA}Q_1 \equiv Q_2 e^{ih^{-1}tA}Q_1 $ modulo operators with $O(h^\infty)$-norms\,\footnote{\label{foot-6} By default, operator norm is from $\sL^2$ to $\sL^2$.} and with $Q_2$ also supported in the small vicinity of $\bar{\xi}$ and equal $1$ in the vicinity of $\supp (Q_1)$.

Then on $\supp(Q_2)$ operator is microhyperbolic with respect to vector $\ell $ and we can employ the proof of
Theorem~\ref{monsterbook-thm-2-1-2}  of \cite{IvrMB} again, this time with $\phi (x,t)= \ell x- \epsilon_0 t$. For further details see Chapter~\ref{monsterbook-sect-4} of \cite{IvrMB}.
\end{proof}

Then in virtue of (\ref{eqn-1.8}) with $t=T_*$ (which is also due to the microhyperbolicity condition) (\ref{eqn-1.8}) also holds with $T=T^*$  and applying H\"ormander's Tauberian theorem we arrive to the remainder estimate $Ch^{1-d} T^{*\,-1}= C\varepsilon h^{1-d}$, thus proving the following theorem:

\begin{theorem}\label{thm-2.4}
Let $A=A^0(hD)+\varepsilon B(x,hD)$ with $A^0$ satisfying conditions \textup{(\ref{eqn-1.1})}, \textup{(\ref{eqn-1.2})} and \textup{(\ref{eqn-1.4})} and $B$ satisfying conditions \textup{(\ref{eqn-1.1})}. Then
\begin{equation}
e_h(x,x,\tau) =\sum_{0\le n\le M}\kappa_n(x,\tau)h^{-d+n}+O(\varepsilon h^{1-d})
\label{eqn-2.9}
\end{equation}
provided $\varepsilon \ge h^{M}$, $|\tau-\lambda|\le \epsilon$.
\end{theorem}

From now on we discuss only Theorem~\ref{thm-1.4}.
\enlargethispage{\baselineskip}

\begin{remark}\label{rem-2.5}
\begin{enumerate}[label=(\roman*), wide, labelindent=0pt]
\item\label{rem-2.5-i}
It suffices to prove asymptotics
\begin{equation}
e_h(x,x,\tau) =\sum_{0\le n\le M}\kappa_n(x,\tau)h^{-d+n}+O(h^{-d+M})
\label{eqn-2.10}
\end{equation}
with arbitrarily large fixed $M$.   To do so we will use the \emph{hyperbolic operator method\/} (which we implement as  \emph{semiclassical Schr\"odinger operator method\/}) with maximal time $T^*=h^{-M}$.

\item\label{rem-2.5-ii}
Then we can replace operator $B$ by operator $B'$, provided operator norm of $B-B'$ from $\sH^{m}$ to $\sL^2$ does not exceed $Ch^{3M}$.

Indeed, let $A'=A^0 +\varepsilon B'$.  Due to Remark~\ref{rem-2.1} we need to compare only $Q_1 e^{ih^{-1}tA'}Q_1$ and $Q_1 e^{ih^{-1}tA}Q_1$. Observe that  due to (\ref{eqn-1.2})
\begin{gather*}
\3 e^{ih^{-1}tA}Q_1- Q_2 e^{ih^{-1}tA}Q_1\3_k \le C_{k,s} h^s
\end{gather*}
with arbitrarily large $k,s$, where $\3 .\3_k$ denotes an operator norms  from $\sL^2$ to $\sH^k$ provided
$\supp(Q_j)\subset \{\xi\colon A^0(\xi)\le 2jc\}$ and $Q_2=1$ in $\{\xi\colon A^0(\xi)<3c\}$. The same is true for $A'$ as well.

Then equality
\begin{equation*}
e^{ih^{-1}t A'}- e^{ih^{-1}tA}= ih^{-1}\int_0^t e^{ih^{-1} (t-t')A}(A'-A)e^{ih^{-1}t'A'}\,dt'
\end{equation*}
and restriction $|t|\le T^*$ imply that
$\3 \bigl(e^{ih^{-1}t A'}- e^{ih^{-1}tA}\bigr) Q\3_k$
does not exceed $C_{k,s}h^{s} + Ch^{-1-M} \3 Q_2 (B-B') Q_2\3_k$. 

Finally, observe that
$\3 Q_2 (B-B') Q_2\3_k \le C_{k} h^{-k-m } \3 (B-B')\3'_k$ where $\3.\3'_k$  denotes an operator norm from $\sH^m$ to $\sL^2$.

\item\label{rem-2.5-iii}
Since $\N_h(\tau )$ could be defined equivalently as
\begin{equation}
\N_h(\lambda)=  \lim _{\ell \to \infty} (\mes (\ell X))^{-1} \N_{h} (\lambda, \ell X)
e(x,x,\lambda)\,dx,
\label{eqn-2.11}
\end{equation}
where $\N_{h} (\lambda,  X)$ is an eigenvalue counting function for operator $A$ in $X$ with the Dirichlet (or Neumann--does not  matter) boundary conditions on $\partial X$, for $\N_h(\tau )$ we can arrive to the same conclusion from the variational arguments.

\item\label{rem-2.5-iv}
First such replacement will be $B'\coloneqq B'_{(L,\omega)}$ from Condition \ref{cond-B} with
$\omega= h^{-\sigma}$, arbitrarily small $\sigma>0$ and  $L=3M/\sigma$.

\emph{So, from now   $\Theta$ and $B$ are effectively replaced by $\Theta'\coloneqq \Theta '_{(L,\omega)}$ and $B'_{(L,\omega)}$ correspondingly\/}. 
\end{enumerate}
\end{remark}

\section{Gauge transformation}
\label{sect-2.2}

Consider now the ``gauge'' transformation $A\mapsto  e^{-i\varepsilon h^{-1} P} A e^{i\varepsilon h^{-1} P}$ with $h$-pseudodifferential operator $P$. Observe that
\begin{multline}
e^{-i\varepsilon h^{-1} P} A e^{i\varepsilon h^{-1} P} = A-i \varepsilon h^{-1}[P,A]
+\sum_{2\le n \le K-1} \frac{1}{n!}(-i\varepsilon  h^{-1})^n \Ad^n_P (A)\\
+
 \int_0^1 \frac{1}{(K-1)!} (1-s)^{K-1}(-i\varepsilon  h^{-1})^K  e^{-i\varepsilon h^{-1} sP}\Ad^K_P(A) e^{i\varepsilon h^{-1} sP}\,ds,
\label{eqn-2.12}
\end{multline}
where $\Ad^0 _P(A)=A$ and $\Ad^{n+1}_P(A)=[P,\Ad^n_P(A)]$ for $n=0,1,\ldots$.

Thus \emph{formally\/} we can compensate $\varepsilon B$, taking
\begin{gather}
P =\sum_{\theta} i h \bigl( A^0(\xi+\theta h/2)- A^0(\xi-\theta h/2)  \bigr) ^{-1} b_\theta (\xi)e^{i\langle \theta,x\rangle},
\label{eqn-2.13}\\
\shortintertext{so that}
ih^{-1}[P, A^0] = B \implies ih^{-1}[P,A]= B + i\varepsilon h^{-1}[P,B].
\label{eqn-2.14}
\end{gather}

Then perturbation $\varepsilon B$ is replaced by  $\varepsilon ^2B'$, which is the  right hand expression in (\ref{eqn-2.12}) minus $A^0$, i.e.
\begin{equation}
 B'=-i  h^{-1}[P,B]
+\sum_{2\le n \le K-1} \frac{1}{n!}\varepsilon^{n-2}(-i  h^{-1})^n \Ad^n_P (A),
 \label{eqn-2.15}
 \end{equation}
where we ignored the remainder.

New perturbation, again formally, has a magnitude of $\varepsilon^2$. Repeating this process we will make a perturbation negligible.

\begin{remark}\label{rem-2.6}
However, we need to address the following issues issues:
\begin{enumerate}[label=(\roman*), wide, labelindent=0pt]

\item\label{rem-2.6-i}
Denominator $h^{-1}\bigl( A^0(\xi+\theta h/2)- A^0(\xi-\theta h/2)  \bigr)= \langle \nabla _\xi A^0,\theta \rangle +O(h^{1-\sigma}) $  could be small.

\item\label{rem-2.6-ii}
In $B'$ set $\Theta '$ increases: $\varepsilon ^2B'=\varepsilon ^2 B'_2 + \varepsilon ^3 B'_3 +\ldots + \varepsilon^M B'_M$, where for $B'_j$ the frequency set is $\Theta'_j$ (the arithmetic sum of $j$ copies of $\Theta'$).

\item\label{rem-2.6-iii}
We need to prove that the remainder is negligible.

\item\label{rem-2.6-iv}
This transformation was used in Section 9 of \cite{ParSht3} (etc); in contrast to these papers we use Weyl quantization instead of $pq$-quantization, and have  therefore
$\bigl( A^0(\xi+\theta h/2)- A^0(\xi-\theta h/2)\bigr)$  instead of
$\bigl( A^0(\xi+\theta h)- A^0(\xi) \bigr)$.
\end{enumerate}
\end{remark}

\section{Non-resonant zone}
\label{sect-2.3}

\subsection{Gauge transformation}
\label{sect-2.3.1}

One can see easily that if inequality
\begin{equation}
|\langle \nabla _\xi A^0(\xi),\theta \rangle| \ge \gamma \coloneqq \varepsilon^{\frac{1}{2}}h^{-\delta} \label{eqn-2.16}
\end{equation}
holds  for all $\theta \in \Theta'_K$, then the terms could be estimated by $h^{\delta n}$ and our construction works with $K=3M/\delta$. Here and below without any loss of the generality we assume that $\varepsilon \ge h$; so, in fact,
\begin{equation}
h^\vartheta \ge \varepsilon \ge h.
\label{eqn-2.17}
\end{equation}
Indeed, if $P=P(x,hD)$ has the symbol, satisfying 
\begin{gather}
|D_\xi ^\alpha D_x^\beta P |\le c_{\alpha\beta} \gamma^{-1-|\alpha|}
\qquad \forall \alpha, \beta,\label{eqn-2.18}\\
\intertext{then $B'=\varepsilon h^{-1} [P,B]$ has a symbol, satisfying}
|D_\xi ^\alpha D_x^\beta B' |\le c'_{\alpha\beta} \varepsilon \gamma^{-2-|\alpha|}
\qquad \forall \alpha, \beta,
\label{eqn-2.19}
\end{gather}
so indeed $\varepsilon'= \varepsilon^2 \gamma^{-2}$.

Then we can eliminate a perturbation completely, save terms with the frequency $0$, both old and new. The set of $\xi$ satisfying (\ref{eqn-2.16}) for all $\theta \in \Theta'_K$ we call \emph{non-resonant zone\/} and denote by $\cZ$. Thus, we arrive to

\begin{proposition}\label{prop-2.7}
Let $Q=Q(hD)$ with the symbol supported in $\cZ\cap \Omega$ and satisfying \textup{(\ref{eqn-2.3})} 

Then there exists a pseudo-differential operator $P=P(x,hD)$ with the symbol, satisfying \textup{(\ref{eqn-2.18})} and such that
\begin{gather}
\bigl(e^{-i\varepsilon h^{-1}P} A e^{i\varepsilon h^{-1}P} -A'' \bigr) Q \equiv 0
\label{eqn-2.20} \\
\shortintertext{with}
A'' = A^0(hD) + \varepsilon B''_0(hD)
\label{eqn-2.21}
\end{gather}
modulo operator from $\sH^m$ to $\sL^2$ with the operator norm $O(h^{3M})$.
\end{proposition}

\begin{remark}\label{rem-2.8}
\begin{enumerate}[label=(\roman*), wide, labelindent=0pt]

\item\label{rem-2.8-i}
This proposition is similar to Lemma 9.3 of \cite{ParSht3}. However, in contrast to \cite{ParSht1, ParSht2, ParSht3, MorParSht}, after it is proven we do not write asymptotic decomposition there, but simply prove that singularities do not propagate with respect to $\xi$ there.

\item\label{rem-2.8-ii}
It is our second replacement of operator $A$; recall that the first one was based on Condition~\ref{cond-B}, and now we ignore the remainder after transformation, which is justified by Remark~\ref{rem-2.5}\ref{rem-2.5-i}.
\end{enumerate}
\end{remark}

\subsection{Propagation}
\label{sect-2.3.12}

\begin{proposition}\label{prop-2.9}
Let $Q_j=Q_j(hD)$ with the symbols,  satisfying \textup{(\ref{eqn-2.3})} and let symbol of $Q_1$ be supported in $\cZ\cap \Omega$.

Let   $\dist(\supp(Q_1),\,\supp(Q_2))\ge c\gamma$. Then
\begin{equation}
\|Q_2e^{ih^{-1}tA}Q_1\|=O(h^{2M}) \qquad \text{as\ \ } |t|\le T^*=h^{-M}.
\label{eqn-2.22}
\end{equation}
\end{proposition}

\begin{proof}
One can prove easily that the operator norms of  $Q_2e^{ih^{-1}tA ''}Q_1$ and  $Q_2e^{\pm i\varepsilon h^{-1}P }Q_1$ are $O(h^{2M})$. We leave all easy details to the reader.
\end{proof}

\section{Resonant zone}
\label{sect-2.4}

Consider now \emph{resonant zone}
\begin{gather}
\Lambda  \coloneqq \bigcup_{\theta \in \Theta'_K\setminus 0} \Lambda (\theta),
 \label{eqn-2.23}\\
\intertext{where $\Lambda (\theta)$ is the set of $\xi $, violating (\ref{eqn-2.16})  for given $\theta $:}
\Lambda (\theta)=\Lambda_\delta(\theta)\coloneqq
\{\xi\colon   |\langle \nabla _\xi A^0(\xi),\theta \rangle| \ge \gamma=c \varepsilon^{\frac{1}{2}}h^{-\delta}\}.
 \label{eqn-2.24}
\end{gather}

\subsection{Case $d=2$}
\label{sect-2.4.1}

We start from the easiest case $d=2$ (in the trivial case $d=1$ there is no resonant zone).  Observe that  due to assumption (\ref{eqn-1.16}) for each $\theta$
\begin{equation}
\mes_1(\Lambda(\theta)\cap \Sigma_\lambda)\le C \gamma .
\label{eqn-2.25}
\end{equation}

Further,  $\# \Theta'_K \le Ch^{-\sigma }$ (as $h\le h_0(K,\sigma)$) due to   Condition~\ref{cond-C}.   Thus
$\mes_1(\Lambda \cap \Sigma_\lambda)\le \gamma h^{ -\sigma}$. Recall, that   $\sigma>0$  is arbitrarily small.

Since due to  Proposition~\ref{prop-2.9}, the propagation which starts in the non-resonant zone $\cZ$ remains there\footnote{\label{foot-7} May be, with different constant $c$ in the definition of $\gamma$.} we conclude that the  propagation which is started in some connected component of the resonant zone also remains there\footref{foot-7}.

Thus,  $\nabla _\xi A^0(\xi)$ does not change by more than $\gamma h^{-\sigma }$ and since $\sigma $ ais arbitrarily small we conclude that (\ref{eqn-2.22}) also holds for $Q_1$, supported in the resonant zone. Therefore

\begin{claim}\label{eqn-2.26}
Estimate (\ref{eqn-2.22}) holds for all $Q_1$, $Q_2$ satisfying (\ref{eqn-2.3}) and
\begin{equation}
\dist (\supp(Q_1),\,\supp (Q_2))\ge \gamma.
\label{eqn-2.27}
\end{equation}
\end{claim}

\begin{remark}\label{rem-2.10}
\begin{enumerate}[label=(\roman*), wide, labelindent=0pt]
\item\label{rem-2.10-i}
In the proof of Theorem~\ref{thm-1.4} we need only to have estimate (\ref{eqn-2.22}) holding for all $Q_1$, $Q_2$ satisfying (\ref{eqn-2.3}) and (\ref{eqn-2.27}) with arbitrarily small constant $\gamma$.

\item\label{rem-2.10-ii}
Then for $d=2$ we can replace assumption (\ref{eqn-1.16}) by
\begin{claim}\label{eqn-2.28}
$\varkappa (s)$ (a curvature of $\Sigma_\lambda$, naturally parametrized by $s$) has zeroes only of the finite order.
\end{claim}
Indeed, then (\ref{eqn-2.25}) will be replaced by $\mes_1(\Lambda(\theta)\cap \Sigma_\lambda)\le C \gamma ^{\nu}$,
$\nu = 1/(q+1)$ with $q$ the maximal order of zeroes of $\varkappa(s)$.
\end{enumerate}
\end{remark}

\subsection{General case: gauge transform}
\label{sect-2.4.2}

Consider now the general case $d\ge 2$.   In this case due Conditions~\ref{cond-A}, \ref{cond-C} and \ref{cond-D} we can cover
$\Lambda\cap \Omega_\tau$ by $\Lambda^*$,
\begin{equation}
\Lambda\cap \Omega_\tau\subset \Lambda^*=\bigcup_{1\le j\le d-1} \Lambda^*_{j} ,
 \label{eqn-2.29}
\end{equation}
defined as:
\begin{claim} \label{eqn-2.30}
Let $\xi \in \Omega_\tau$; then $\xi \in \Lambda_j^*$ iff there exist $\theta_1,\ldots ,\theta_j\in \Theta'_K$ which are linearly independent and such that $\xi \in \Lambda _{\delta_j} (\theta_k)$ for all $k=1,\ldots, j$,
\end{claim}
where $0<\delta=\delta_1 <\delta_2<\ldots <\delta_{d-1}$ are arbitrarily fixed and we chose sufficiently small
$\sigma>0$ afterwards.

Further, due to Conditions~\ref{cond-A}, \ref{cond-C}, \ref{cond-D}  and (\ref{eqn-1.16})
$\Lambda^*_{d-1}\cap \Omega_\tau$ could be covered by no more than $\gamma_{d-1}$-vicinities of some points
$\xi_\iota$,  $\iota=1,\ldots, \omega^{g}$, $g=g(d)$. Recall that
$\Omega_\tau\coloneqq \{\xi\colon | A^0(\xi)-\tau |\le C_0\varepsilon + h^{1-\varsigma}\}$.

Consider some connected component $\Xi$ of $\Lambda^*_j$. Let some point $\bar{\xi}$ of it belong to
$\bigcap_{1\le k \le j} \Lambda_{\delta_j} (\theta_k)\cap \Omega_\tau$ with linearly independent
$\theta_1,\ldots,\theta_j$. Observe that
$\diam (\bigcap_{1\le k \le j} \Lambda_{\delta_j} (\theta_k)\cap \Omega ) \le c\gamma_j$ due to strong convexity assumption (\ref{eqn-1.16}). Then this set either intersects or does not intersect with $\Lambda^*_{j+1}\cap \Omega$. In the former case we include it to $\Lambda_{j+1}^*$ and exclude it from $\Lambda_j^*$.

After we redefined  $\Lambda_j^*$ we arrive to the following proposition:

\begin{proposition} \label{prop-2.11}
Equation \textup{(\ref{eqn-2.29})} still holds where now each connected component $\Xi $ of $\Lambda_j^*$ has the following properties:

\begin{enumerate}[label=(\roman*), wide, labelindent=0pt]

\item\label{prop-2.11-i}
$\diam \Xi \le c \gamma_j$.

\item\label{prop-2.11-ii}
There exist linearly independent $\theta_1,\ldots,\theta_j\in \Theta'_K$, such that for each $\xi\in \Xi $
$|\langle \nabla_\xi A^0(\xi) ,\theta\rangle|\le c_j \gamma_j $
for  all $\theta \in \fV \cap (\Theta'_K\setminus 0)$ and
$|\langle \nabla_\xi A^0(\xi) ,\theta\rangle|\ge \epsilon_j \gamma_{j+1}$
for all $\theta \in \Theta'_K\setminus \fV)$  with $\fV=\Span (\theta_1,\ldots,\theta_j)$.
\end{enumerate}
\end{proposition}

Now we generalize Proposition~\ref{prop-2.7}:

\begin{proposition}\label{prop-2.12}
Let $Q=Q(hD)$ with the symbol supported in the connected component $\Xi $ of $\Lambda_j^*$, corresponding to
subspace $\fV$, and satisfying \textup{(\ref{eqn-2.3})}. Then there exists a pseudo-differential operator $P=P(x,hD)$ with the symbol, satisfying \textup{(\ref{eqn-2.18})} and such that
\begin{gather}
\bigl(e^{-i\varepsilon h^{-1}P} A e^{i\varepsilon h^{-1}P} -A'' \bigr) Q \equiv 0
\label{eqn-2.31}
\end{gather}
modulo operator from $\sH^m$ to $\sL^2$ with the operator norm $O(h^{3M})$, where
$A'' =A^0 +\varepsilon B'' (x,hD)$, where $B''$ is an operator with Weyl symbol
\begin{gather}
B'' (x,\xi)= \sum_{\theta\in \Theta '_K\cap \fV} b_{\fV, \theta}(\xi) e^{i\langle \theta,x\rangle}.
\label{eqn-2.32}
\end{gather}
\end{proposition}

\begin{proof}
The proof obviously generalizes the proof of Proposition~\ref{prop-2.7}. We eliminate all $\theta \notin \fV$ exactly in the same way as it was done there.
\end{proof}

\subsection{General case: propagation}
\label{sect-2.4.3}

\begin{proposition}\label{prop-2.13}
Let $Q_j=Q_j(hD)$ with the symbols,  satisfying \textup{(\ref{eqn-2.3})} and let symbol of $Q_1$ be supported in
$\Lambda^*_j$ .

Let  $\dist(\supp(Q_1),\,\supp(Q_2))\ge C_0\gamma_j$. Then $\|Q_2e^{ih^{-1}tA}Q_1\|=O(h^{2M})$  for
$|t|\le T_*=h^{-M}$.
\end{proposition}

\begin{proof}
In virtue of Proposition~\ref{prop-2.9} it is sufficient to consider $\supp (Q_1)$ belonging to the connected component $\Xi'$ of $\Lambda_j^*$. Indeed, the values of $\delta_1,\ldots, \delta_{d-1}$ are arbitrarily small.

One can prove easily that the operator norm of    $Q_2e^{\pm i\varepsilon h^{-1}P }Q_1$ are $O(h^{2M})$.
We need to prove that the operator norm of    $Q_2e^{\pm i h^{-1} t A'' }Q_1$ is also $O(h^{2M})$. In the coordinates
$(x';x'')\in \fV \oplus (\bR^d\ominus \fV)$ we observe that the propagation speed is  only along $\fV$ as long as it remains in $\epsilon \gamma_j$ vicinity of $\supp (Q_1)$. The proof is similar to the proof of Proposition~\ref{prop-2.2} and we leave it to the reader.

However propagation is confined to $\Omega '_\tau\coloneqq \{\xi\colon |A^0(\xi)-\tau|\le C\varepsilon + 2h^{1-\varsigma}\}$) and due to (\ref{eqn-1.16}) it remains in that vicinity as $\varsigma <\delta$.
\end{proof}

Now we arrive to the following proposition:

\begin{proposition}\label{prop-2.14}
Let $Q_1, Q_2$ satisfy \textup{(\ref{eqn-2.3})} and $\supp (Q_1) \subset \Omega$. Then for $T_*\le T\le T^*$
\begin{equation}
 F_{t\to h^{-1}\tau}\bigl(\chi _T(t)Q_{2x} u(x,y,t)\,^t\!Q_{1y}\bigr) = O(h^{2M}).
 \label{eqn-2.33}
\end{equation}
\end{proposition}

\begin{proof}
It is standard, due to Proposition~\ref{prop-2.13}, microhyperbolicity condition and the results of Chapter~\ref{sect-2} of \cite{IvrMB} we conclude that if $|\ell|=1$ and
\begin{gather}
\langle \ell, \nabla_\xi A^0(\xi) \rangle \ge \epsilon _0 \qquad \forall \xi \in \supp(Q_1)
\label{eqn-2.34}\\
\shortintertext{and}
\langle \ell, x-y\rangle \le \epsilon_1 T \qquad \forall x \in \supp(\phi_1), \ y\in \supp (\phi_2),
\label{eqn-2.35}
\end{gather}
then $\| \phi_2 e^{i h^{-1}t A}Q_1 \phi_2\|=O(h^{2M})$ for $T\le t\le 2T$.

This implies (\ref{eqn-2.33}) provided  $\diam(\supp (Q_1))\le \epsilon $. But then for (\ref{eqn-2.33}) we can drop this assumption.
\end{proof}

\section{End of the proof}
\label{sect-2.5}
Now we conclude that
\begin{equation}
 F_{t\to h^{-1}\tau}\bigl([\bar{\chi} _{T}(t)-\bar{\chi} _{T_*}(t)]Q_{2x} u(x,y,t)\,^t\!Q_{1y}\bigr)\bigr|_{x=y}= O(h^{2M})
 \label{eqn-2.36}
 \end{equation}
 and since
 \begin{multline}
F_{t\to h^{-1}\tau}\bigl(\bar{\chi} _{T}(t) Q_{2x} u(x,y,t)\,^t\!Q_{1y}\bigr)\bigr|_{x=y}= \\
\sum_{0\le n\le M} \kappa' _n (x,\varepsilon) h^{1-d+n} + O(h^{M+1})
 \label{eqn-2.37}
\end{multline}
holds for $T=T_*$, it also holds for $T=T^*$.

Finally, H\"ormander's Tauberian theorem implies Theorem~\ref{thm-1.4}.

\newpage

\chapter{Generalizations and Discussion}
\label{sect-3}
\section{Matrix operators}
\label{sect-3.1}

Consider now $n\times n$-matrix operators $A^0$ and $B$; then (\ref{eqn-1.2}) should be understood in the matrix sense. Assume that

\begin{claim}\label{eqn-3.1}
Symbol $A^0(\xi)$ has only simple eigenvalues $a^0_1(\xi),\ldots ,a^0_n(\xi)$, which also satisfy (\ref{eqn-1.4}) and (\ref{eqn-1.16}).
\end{claim}

Then there exists a unitary transformation $R^0=R(\xi)$, such that
$R^{0\,\dag}(\xi)A^0(\xi)R^0(\xi)=\diag (a^0_1(\xi),\ldots , a^0_n(\xi))$.

Then one can prove easily, that there exists a unitary operator $R(x,hD)=R^0(hD)+\varepsilon R'(x,D)$, such that
$R^{*}AR =\diag (a_1,\ldots , a_n)$, where $a_j=a_j(x,hD)=a_j^0(hD)+ \varepsilon b_j (x,hD)$ (and we assume as before that (\ref{eqn-2.17}) holds.

If Conditions \ref{cond-A}--\ref{cond-D} are fulfilled for $A(x,hD)$, then they are also fulfilled for $a_j(x,hD)$ and we can apply the same propagation arguments as before and Theorem \ref{thm-1.4} extends to such operators provided conditions (\ref{eqn-1.4}) and (\ref{eqn-1.16}) are fulfilled for $a_j(x,hD)$ with $j=1,\ldots,n$.

Let us replace (\ref{eqn-1.2}) by more general ellipticity assumption
\begin{equation}
|A^0(\xi)v|\ge \epsilon |\xi|^m |v|\qquad \forall v\in \bC^n\ \forall \xi\colon |\xi|\ge C_0.
\label{eqn-3.2}
\end{equation}

Then we cannot restrict $e(x,y,\lambda)$ to $x=y$ but we can restrict $e(x,y,\lambda,\lambda')$, the Schwartz kernel of the difference of the corresponding projectors.

Theorem~\ref{thm-1.4} trivially extends to such operators, if instead of $e(x,x,\lambda)$ we consider
$e(x,x,\lambda,\lambda')$ provided conditions (\ref{eqn-1.4}) and (\ref{eqn-1.16}) are  fulfilled for $a_j(x,hD)$
with $j=1,\ldots,n$ and for both $\lambda $ and $\lambda'$. It also extends to
\begin{equation}
\int e(x,y,\lambda,\lambda')\phi (\lambda') \,d\lambda',\qquad  \phi\in \sC_0^\infty (\bR),
\label{eqn-3.3}
\end{equation}
provided conditions (\ref{eqn-1.4}) and (\ref{eqn-1.16}) are fulfilled for $a_j(x,hD)$ with $j=1,\ldots, n$ for $\lambda $.

\begin{remark}\label{rem-3.1}
Our reduction construction fails in the case of a scalar operator $A^0$ and a matrix operator $B$ unless either $\varepsilon =h^{1+\delta}$ or the principal symbol of $B$ satisfies some very restrictive condition.
Therefore for a matrix operator $A^0$ with the eigenvalues of $A^0(\xi)$ of constant multiplicities our construction works only under similar assumptions.
\end{remark}

\section{Perturbations}
\label{sect-3.2}

Consider operators in question, perturbed by $\varepsilon V(x,hD)$ where $V(x,\xi)$ decays as $|x|\to \infty$. Such perturbations do not affect $\N_h(\lambda)$, but they do affect $e_h(x,x,\lambda)$.
\subsection{Decaying perturbations}
\label{sect-3.2.1}

 We start from the easy case
\begin{gather}
A= A^0(hD) + \varepsilon V(x,HD),
\label{eqn-3.4}\\
\shortintertext{where}
|D^\alpha _\xi D_x^\beta V(x,\xi)|\le  c_{\alpha\beta}(|\xi|+1)^m (|x|+1)^{-\delta-|\beta|} \qquad \forall \alpha, \beta\ \forall x, \xi.
\label{eqn-3.5}
\end{gather}

First of all, we claim that
\begin{claim}\label{eqn-3.6}
Under assumption (\ref{eqn-3.7}) below the propagation speed with respect to $\xi$ does not exceed $c\varepsilon (|x|+1)^{-\delta}$.
\end{claim}

Indeed, note first that due to  Proposition~\ref{prop-2.2} the propagation speed with respect to $\xi$ does not exceed $c\varepsilon$.  Next, consider domain $\{x\colon |x|\asymp r\}$ with $r\ge 1$. Scaling $x\mapsto x/r$, $t\mapsto t/r$ we get a domain  $\{x\colon |x|\asymp 1\}$, $h\mapsto \hbar =h/r$  and we need to prove that after this  scaling the propagation speed with respect to  $\xi$ does not exceed $\nu = c\varepsilon r^{-\delta}$, on the time interval $\{t\colon |t| \le 1\}$.

To prove this we can apply Proposition~\ref{prop-2.2} but ewe need to have the microlocal uncertainty principle fulfilled: $\nu \ge \hbar ^{1-\sigma}$ with $\sigma >0$, where $\nu$ is a shift with respect to $\xi$. This inequality is equivalent to  $\varepsilon r^{-\delta}\ge h^{1-\sigma} r^{-1+\sigma}$ i.e. 
$\varepsilon r^{1-\sigma-\delta}\ge h^{1-\sigma}$ and it suffice to have
\begin{equation}
\delta <1, \qquad \varepsilon \ge h^{1-\sigma} \quad \text{with\ \ } \sigma>0.
\label{eqn-3.7}
\end{equation}

Consider now $\xi $ in the vicinity of $\bar{\xi}$ and $x$  with $|x|\le c$. Then as long as $|\xi-\bar{\xi}|\le \epsilon$ with small enough constant $\epsilon >0$, evolution goes away from $0$ with the speed $\asymp 1$, so we are in the zone $\{x\colon |x|  \asymp |t|\}$ and in this zone the propagation speed with respect to $\xi$ does not exceed $c \varepsilon r^{-1-\delta}$, and therefore $|\xi-\bar{\xi}| \le c\varepsilon \int_1^\infty t^{-1-\delta}\,dt \le  c\varepsilon$ and this is less that $\epsilon/2$ as $\varepsilon \le \epsilon_0$.

We can also consider evolution which starts from $x$ with $|x|\ge 1$. Then the same arguments work albeit with $r \asymp |t-t^*|$ for some $t^*$ with $|t^*|\le c|x|$.

Then we arrive to

\begin{theorem}\label{thm-3.2}
Consider operator~\textup{(\ref{eqn-3.4})} with $V$ satisfying~\textup{(\ref{eqn-3.5})}. Let microhyperbolicity condition \textup{(\ref{eqn-1.4})} on the energy level $\lambda$ be fulfilled and $\varepsilon \le \epsilon_0$. Then  the complete spectral asymptotics \textup{(\ref{eqn-1.29})} holds.\end{theorem}


\subsection{Hybrid perturbations}
\label{sect-3.2.2}

Now we consider the hybrid operators, containing both $\varepsilon B$ and $\varepsilon V$. However, trying to eliminate  $\varepsilon B$  by the same approach as in Subsubsection~\ref{sect-2.4.2}.2, we get an another type of terms, and it is only natural to consider them being in the operator from the beginning:
\begin{gather}
A = A^0 (hD) + \varepsilon \bigl( B(x,hD) + V(x,hD)\bigr),
\label{eqn-3.8}\\
\shortintertext{where}
V(x,\xi)=\sum _{\theta \in \Theta} e^{i\langle \theta, x\rangle} V_\theta (x,\xi),
\label{eqn-3.9}\\
|D^\alpha _\xi D_x^\beta V(x,\xi)|\le  c_{\alpha\beta}(|\xi|+1)^m (|x|+1)^{-\delta} \quad \forall \alpha, \beta\ \forall x, \xi.
\label{eqn-3.10}
\end{gather}

We impose condition
\begin{condition}\label{cond-E}
For each $\omega$ and $L$ for the same set $\Theta'$ as before there exists
\begin{gather}
V'(x,\xi)=\sum _{\theta \in \Theta'} e^{i\langle \theta, x\rangle} V'_\theta (x,\xi),
\label{eqn-3.11}\\
\shortintertext{such that}
\| D_x^\alpha D_\xi^\beta \bigl(V -V' \bigr) \|_{\sL^\infty}\le   \omega^{-L}(|\xi|+1)^m\label{eqn-3.12}\\
\shortintertext{and}
| D_x^\alpha D_\xi^\beta V'_\theta  |\le  c_{Ls\alpha\beta}(|x|+1)^{-1-\delta -|\alpha|}(|\theta|+1)^{-s}
\label{eqn-3.13}
\\[5pt]
\hphantom{\hskip250pt}\forall \alpha, \beta\colon |\alpha| \le L, |\beta|\le L \ \forall s.
\notag
\end{gather}
\end{condition}

\subsection{Non-resonant zone}
\label{sect-3.2.3}
We deal with the purely exponential terms in our standard way and with the hybrid terms as if they were purely exponential (i.e. as if $V'_\theta$ were not depending on $x$), then a new kind of terms will be produced: they acquire factor
$h(A^0(\xi +\theta h/2)-A^0(\xi -\theta h/2))^{-1}$ and the derivative with respect to $x$ to $V'_\theta$.

Eventually we end up with the operator of the same type (\ref{eqn-3.8}) with $B(x,\xi)$ replaced by
$B''(\xi)$ and with $V_\theta(x,\xi)$ replaced by  $V''_\theta (x,\xi)$, such that
\begin{equation*}
|D^\alpha_\xi D^\beta_x V''_\theta (x,\xi)|\le C_{n\alpha\beta}\varepsilon^{k+1}\gamma^{-2k-n-|\alpha|} (|x|+1)^{-n-\delta-|\beta|} 
\end{equation*}
with $n+k\ge 3K$.

Then
\begin{equation*}
|D^\alpha_\xi D^\beta_x \bigl[V''_\theta (x,\xi)e^{i\langle \theta, x\rangle}\bigr]|\le C_{s\alpha\beta}\varepsilon^{k+1}\gamma^{-2k-n-|\alpha|} (|x|+1)^{-n-\delta} (|\theta|+1)^{|\beta|};
\end{equation*}
recall that $|\theta|\le CKh^{-\sigma}$.

Let us pick up $\gamma =h^{\delta}$ with $\delta =\vartheta /6K$. Then,  ignoring therms with $k\ge K$ which are negligible, and following the proof of (\ref{eqn-3.6}), we can recover the same statement for the operator after transform, and, finally, to the analogue of Proposition~\ref{prop-2.9}.

\subsection{Resonant zone}
\label{sect-3.2.4}
If $d=2$ we arrive to the analogue of Proposition~\ref{prop-2.2} in the virtue of the we arguments as in Subsubsection~\ref{sect-2.4.1}.1.

If $d\ge 3$ we apply the reduction, similar to one, used in Subsubsection~\ref{sect-2.4.3}.3, and arrive again to operator of the type (\ref{eqn-3.8}) with $B$ replaced by $B''(x,\xi')$ and with $V_\theta(x,\xi)$ replaced by
 $V''_\theta(x,\xi)$. 
 
Then we observe that the shift in direction $\bR^d\ominus\fV$ does not exceed 
$c\varepsilon ^{\delta/2}$ and if it is $\ll \gamma^2$ we arrive to the analogue of Proposition~\ref{prop-2.13}.
It is doable by the choice of really small $\sigma_1<\ldots<\sigma_{d-1}$. Then we arrive to the analogue of Proposition~\ref{prop-2.14} and, finally, to
 
\begin{theorem}\label{thm-3.3}
Let $A$ be a self-adjoint operator \textup{(\ref{eqn-3.8})},  where $A^0$ satisfies
\textup{(\ref{eqn-1.1})}, \textup{(\ref{eqn-1.2})}, \textup{(\ref{eqn-1.4})} and \textup{(\ref{eqn-1.16})} and  $B$ satisfies \textup{(\ref{eqn-1.1})}, $V$ satisfies \textup{(\ref{eqn-3.9})} and \textup{(\ref{eqn-3.10})}.

Let Conditions \ref{cond-A}--\ref{cond-E} be fulfilled. Then for $|\tau-\lambda|<\epsilon$, $\varepsilon \le h^{\vartheta}$, $\vartheta  >0$ asymptotics \textup{(\ref{eqn-1.29})} holds.
\end{theorem}

\section{Differentiability}
\label{sect-3.3}

It also follows from Corollary~\ref{cor-1.5} that
\begin{equation}
\frac{1}{\nu}\Bigl[\N_{h,\varepsilon}(\tau + \nu) -\N_{h,\varepsilon}(\tau)\Bigr] =
\frac{1}{\nu}\Bigl[\cN_{h,\varepsilon}(\tau + \nu) -\cN_{h,\varepsilon}(\tau)\Bigr]+O(h^\infty)
\label{eqn-3.14}
\end{equation}
provided $\nu \ge h^M$, where $\cN_{h,\varepsilon}(\tau)$ is the right-hand expression of (\ref{eqn-1.30}).

The question remains, if (\ref{eqn-3.14}) holds for smaller $\nu$, in particular, if it holds in $\nu\to 0$ limit? If the latter holds, then
\begin{equation}
\frac{\partial\ }{\partial \tau} \N_{h,\varepsilon}(\tau ) =
\frac{\partial\ }{\partial \tau} \cN_{h,\varepsilon}(\tau)+O(h^\infty)
\label{eqn-3.15}
\end{equation}
and we call the left-hand expression the \emph{density of states\/}.

It definitely is not necessarily true, at least in dimension $1$. From now on we consider only asymptotics with respect to $\tau\to +\infty$. Let $A= \Delta + V(x)$ with periodic $V$. It is well-known that for $d=1$ and generic periodic $V$ all spectral gaps are open which contradicts to 
\begin{equation}
\frac{\partial\ }{\partial \tau} \N (\tau ) =
\frac{\partial\ }{\partial \tau} \cN (\tau)+O(\tau^{-\infty}).
\label{eqn-3.16}
\end{equation}

On the other hand, this objection does not work in case $d\ge 2$ since only  several the lowest spectral gaps are open (Bethe-Sommerfeld conjecture, proven in \cite{ParSob}).

Assume for simplicity, that $A=\Delta +V$ has no negative eigenvalues; then we can apply wave operator method\footnote{\label{foot-8} It could be applied without this assumption, but with tweaking.}. We consider $u(x,y,t)$, the Schwartz kernel of $\cos (\sqrt{A} t)$,
\begin{equation}
u(x,y,t) = \int \cos (t\tau)\,d_\tau e(x,y,\tau^2).
\label{eqn-3.17}
\end{equation}
Then, for compactly supported $V$\,\footnote{\label{foot-9} It, probably could be proven for $V$, decaying fast enough at infinity}
\begin{equation}
u(x,y,t) = \left\{\begin{aligned}
&O(e^{-\epsilon |t|}) &&\text{for odd\ \ } d,\\
&O(|t|^{-d}) &&\text{for even\ \ } d
\end{aligned}\right.
\label{eqn-3.18}
\end{equation} 
as $|x|+|y|\le c$, $|t|\to +\infty$ and $\frac{\partial\ }{\partial\tau} e(x,x, \tau^2)$ could be completely restored by inverse $\cos$-Fourier transform, without any Tauberian theorem, and we arrive to asymptotics of $\frac{\partial\ }{\partial\tau} e(x,x, \tau^2)$. Moreover, we can differentiate  
complete asymptotics of the \emph{Birman-Schwinger spectral shift function\/}
\begin{gather}
\xi  (\tau)\coloneqq \int\bigl(e (x,x,\tau^2) - e^0 (x,x,\tau^2)\bigr)\,dx \sim
\sum_{n\ge 0}\bar{\kappa}_n \tau ^{-d+n},
\label{eqn-3.19}\\
\shortintertext{with}
\bar{\kappa}_n\coloneqq\int (\kappa_n(x)-\kappa_n^0)\,dx,
\label{eqn-3.20}
\end{gather}
where $e^0(x,y,\tau)$ and $\kappa_n^0$ correspond to $A^0=\Delta$. In the case of $A=\Delta $ in the exterior of smooth, compact and non-trapping obstacle and $A^0=\Delta$ in $\bR^d$ such asymptotics was derived in \cite{petkov:popov}.

\end{document}